\documentclass{amsproc}
\usepackage{amssymb}

\theoremstyle{plain}

\newtheorem{corollary}{Corollary}

\newtheorem{example}{Example}

\newtheorem{lemma}{Lemma}

\newtheorem{remark}{Remark}

\newtheorem{theorem}{Theorem}
\numberwithin{equation}{section}
\input{tcilatex}

\begin{document}
\title[Completeness and uncertainty principles]{Completeness, special
functions and uncertainty principles over $q$-linear grids}
\author{Lu\'{\i}s Daniel Abreu}
\address{Department of Mathematics, Universidade de Coimbra, Portugal}
\email{daniel@mat.uc.pt}
\thanks{}

\begin{abstract}
We derive completeness criteria for sequences of functions of the form $%
f(x\lambda _{n})$, where $\lambda _{n}$ is the $nth$ zero of a suitably
chosen entire function. Using these criteria, we construct complete
nonorthogonal systems of Fourier-Bessel functions and their $q$-analogues,
as well as other complete sets of $q$-special functions. We discuss
connections with uncertainty principles over $q$-linear grids and the
completeness of certain sets of $q$-Bessel functions is used to prove that,
if a function $f$ and its $q$-Hankel transform both vanish at the points $%
\{q^{-n}\}_{n=1}^{\infty }$, $0<q<1$, then $f$ must vanish on the whole $q$%
-linear grid $\left\{ q^{n}\right\} _{n=-\infty }^{\infty }$.
\end{abstract}

\maketitle

\section{Introduction}

\emph{The Heisenberg uncertainty principle from quantum mechanics.} It can
be reformulated as a proposition saying that if a function is
\textquotedblright small\textquotedblright\ outside an interval of length $T$
and its Fourier transform is \textquotedblright small\textquotedblright\
outside an interval of length $\Omega $, then the product $T\Omega $ must be
bigger than a certain positive quantity. This idea has been used not only in
quantum mechanics, but also in time-frequency analysis, signal recovering
and partial differential equations. Variations of the Heisenberg uncertainty
principle include more quantitative versions and propositions related to the
nature of the support of a function. Integral transforms other than the
Fourier transforms have also been considered and discrete forms of
uncertainty principles constitute a topic of particular interest. Every
uncertainty principle is an instance of a metaproposition which says that a
function and its transform cannot be simultaneously \textquotedblright
small\textquotedblright .

\emph{Completeness of sets of functions.} Over the years, attention has been
given to sequences of functions that, although not being necessarily a basis
of a given space, do however posses the property that every function in that
space can be approximated arbitrarily closely by finite combinations of
those sequences. These sequences are said to be complete and, in the
presence of the classical Fourier setting, they correspond, via Fourier
duality, to uniqueness sets in the Paley-Wiener space.

Following the pioneering work of Paley, Wiener and Levinson, a considerable
amount of research has appeared, concerning completeness properties of the
complex exponentials $\{e^{i\lambda _{n}t}\}$, giving rise to the theory of
nonharmonic Fourier series \cite{Young}, which provide a theoretical
framework for irregular sampling theory.

An important completeness result in the classical Fourier setting states
that, if $\{\lambda _{n}\}$ is the set of zeros of a function of sine type,
then the system $\{e^{i\lambda _{n}t}\}$ is complete in $L^{2}\left[ -\pi
,\pi \right] $ \cite[pag. 145]{Young}\ . In such a context, the set of zeros
of a function of sine type can be seen as a deformation of the set of zeros
of the function sine, $\{\pi n\}$. Using this as a model, it is natural to
try to understand completeness properties of sequences defined in an
analogous way, but replacing the complex exponential with other special
functions. With this in mind, define a sequence $\left\{ f_{n}\right\} $\ of
functions by 
\begin{equation}
f_{n}\left( x\right) =f\left( \lambda _{n}x\right)   \label{main}
\end{equation}%
where $f$ is an entire function and $\lambda _{n}$ is the $nth$ zero of
another entire function $g$. The task is to find conditions in $f$ and $g$
that imply completeness of \ $\{f_{n}\}$. Such a question is particulary
interesting when, for some sequence $\{\lambda _{n}\}$, it is known that the
functions in (\ref{main}) form an orthogonal basis for a given space. In
this case, the functions in the general case can be seen as a deformation of
such a basis. This idea has very classical roots. It is foremost inspired by
Boas and Pollard, who studied in \cite{BP} sequences of nonorthogonal
Fourier-Bessel functions $\{J_{\nu }(\lambda _{n}x)\}$ where $\lambda _{n}$
is not necessarily the $nth$ zero of $J_{\nu }$. A good summary of classical
methods to study general complete systems of special function is Higgins
monograph \cite{Higgins} \ A reading of this monograph and a confrontation
with the revised edition of Young%
\'{}%
s book \cite{Young} gives a historical feeling of how the completeness
problems inspired much of the modern frame theory. The recent developments
concerning expansions in Fourier series on $q$-linear grids \cite{BC}, \cite%
{AM} and the construction of $q$-sampling theorems \cite{A}, \cite{A2}, \cite%
{Ann}, \cite{IZ} motivated the necessity of developing methods to prove
completeness of these systems.

The purpose of our work is twofold. We will first derive completeness
criteria for sequences of the type (\ref{main}) and illustrate them with
several examples involving special and $q$-special functions. As a second
goal of the paper we will obtain uncertainty principles for functions
defined over $q$-linear grids, by proving two statements about a certain $q$%
-analogue of the Hankel transform, introduced by Koornwinder and Swarttouw
in \cite{KoorS}.

\emph{Outline of the  paper:}

Section two recalls some function theoretical definitions and a result is
proved assuring, under certain conditions, the $L^{p}\left[ \mu ,X\right] $, 
$p\geq 1$, completeness of the sequence $\{f\left( \lambda _{n}x\right) \}$,
when $\lambda _{n}$ is the $nth$ zero of $g$, a suitably chosen entire
function, both $f$ and $g$ of order less that one.\ The proof of this
criterion is simple and very classical in nature, using classical entire
function theory. The main argument rests on an application of the Phragm\'{e}%
n-Lindel\"{o}f principle in a proper setting. The general form of the
functions $f$ and $g$ is restricted in a such a way that it fits to many of
the classical and $q$-classical special functions. Within this setting, our
completeness criteria will be illustrated with sets of nonorthogonal
Fourier-Bessel functions, Euler infinite products and Jackson 
\'{}%
s second and third $q$-Bessel functions. The most significative among these
examples is the one using the third Jackson $q$-Bessel function, since the
completeness result can be seen as a deformation of the orthogonal case.

In the third section we obtain a version of the completeness criterion for
functions of order less than two. Specifying the measure $d\mu $ to the
usual $dx$ measure in the real line we apply the result to Bessel functions.
We stress that the completeness criteria from section 2 and 3 are not
contained in Higgins monograph \cite{Higgins} and an extensive search in the
overall available literature seems to sopport that they are new.

In the last section we study uncertainty principles over $q$-linear grids, a
topic that, at a first glance, may seem to be unrelated to the previous one.
However, a relation does exist. The section begins with a brief paragraph
about uncertainty principles. Then, the relevant known facts about the $q$%
-Hankel transform are presented and using a completeness result for $q$%
-Bessel functions of the third type from the previous section, we derive a
vanishing theorem stating that a $L_{q}^{1}(\mathbf{R}^{+})$ function and
its $q$-Hankel transform cannot be both simultaneously supported at the $q$%
-linear grid $\left\{ q^{n}\right\} _{n=-\infty }^{\infty }\cap
(0,1)=\left\{ q^{n}\right\} _{n=1}^{\infty }$, without vanishing on the
equivalent classes of $L_{q}^{1}(\mathbf{R}^{+})$ (that is, on the whole
grid $\left\{ q^{n}\right\} _{n=-\infty }^{\infty }$). The discussion is
then complemented with an uncertainty principle of the type studied in \cite%
{BP} by Donoho and Stark, using the concept of $\epsilon $-concentration.
Here the main instruments used in the proof are a proposition due to de Jeu 
\cite{Jeu} and an estimate on the third Jackson $q$-Bessel function from 
\cite{FHB}. Our uncertainty principles are different from the $q$-analogue
of the Heisenberg uncertainty relations provided in \cite{BFB}.

\section{Completeness criteria}

The unifying theme through this work will be the $L^{p}$-completeness of a
sequence of functions. A sequence of functions $\left\{ f_{n}\right\} $ is
complete in $L^{p}\left[ \mu ,X\right] $ provided the relations 
\begin{equation*}
\int_{X}yf_{n}d\mu =0
\end{equation*}%
for $n=1,2,...$, with $y\in L^{p}\left[ \mu ,X\right] $ and $1/p+1/q=1$,
imply $y=0$ almost everywhere. If $X$ is a finite interval, then $L^{p}\left[
\mu ,X\right] \subset L^{1}\left[ \mu ,X\right] $ , $p\geq 1$, and
completeness in $L^{1}\left[ \mu ,X\right] $ carries with it completeness in 
$L^{p}\left[ \mu ,X\right] $, $p\geq 1$. We will borrow terminology from
Boas and Pollard and say that a set is complete $L\left[ \mu ,X\right] $ if
it is complete in $L^{1}\left[ \mu ,X\right] $.

Some facts from the classical entire function theory will be used in this
section. The maximum modulus of the entire function $f$ is defined as 
\begin{equation*}
M\left( r;f\right) =\max_{\left\vert z\right\vert =r}\left\vert f\left(
z\right) \right\vert
\end{equation*}
and the order of $f$ as 
\begin{equation}
\varrho \left( f\right) =\lim_{r\rightarrow \infty }\frac{\log \log M\left(
r;f\right) }{\log r}.  \label{deford}
\end{equation}
In the case where $f$ is a canonical product with zeros $r_{1},r_{2},...$,
the order of $f$\ is equal to the greatest lower bound of all the $\tau $
for which the series 
\begin{equation*}
\sum\limits_{n=1}^{\infty }\frac{1}{\left\vert r_{n}\right\vert ^{\tau }}
\end{equation*}
converges. From this it is easy to verify that if $A\subset B$ then 
\begin{equation}
\varrho \left[ \prod\limits_{n\in A}^{\infty }\left( 1-\frac{z}{r_{n}}%
\right) \right] \leq \varrho \left[ \prod\limits_{n\in B}^{\infty }\left( 1-%
\frac{z}{r_{n}}\right) \right] .  \label{prop}
\end{equation}
The proof of the main result requires the following form of the Phragm\'{e}%
n-Lindel\"{o}f Principle \cite{Young}\textit{.}

\textbf{\ }\emph{If the order of an entire function }$f$\emph{\ is less than 
}$\sigma $\emph{\ and }$f$\emph{\ is bounded on the limiting rays of an
angle with opening }$\pi /\sigma $\emph{\ then }$f$\emph{\ is bounded on the
region defined by the rays. }

Our general setting is constituted by two nonnegative sequences of real
numbers $(a_{n})$ and $(b_{n})$, defining two entire functions $f$ and $g$
by means of the power series expansions 
\begin{equation}
f(z)=\sum\limits_{n=0}^{\infty }(-1)^{n}a_{n}z^{2n\text{ }}\text{ }
\label{fgseries}
\end{equation}%
and 
\begin{equation}
g(z)=\sum\limits_{n=0}^{\infty }(-1)^{n}b_{n}z^{2n\text{ }}.  \label{gseries}
\end{equation}%
Assume that the zeros of $f$ and $g$ are real, simple, and that there exists
a countable infinite number of them. Denote by $\zeta _{n}$ the $nth$
positive zero of $f$ and denote by $\lambda _{n}$ the $nth$ positive zero of 
$g$. Our first result is the following

\begin{theorem}
\textbf{\ }Let $\mu $\ be a real positive measure. If the order of $f$\ and $%
g$\ are less than one, then the sequence $\{f(\lambda _{n}x)\}$\ is complete 
$L\left[ \mu ,\left( 0,1\right) \right] $\ if, as $n\rightarrow \infty ,$%
\emph{\ } 
\begin{equation}
\frac{a_{n}}{b_{n}}\rightarrow 0\text{ .}  \label{cond}
\end{equation}
\end{theorem}

$\mathbf{Proof.}$ Let $y\in L\left[ \mu ,\left( 0,1\right) \right] $ such
that for $n=1,2,...$%
\begin{equation}
\int_{0}^{1}y\left( x\right) f(\lambda _{n}x)d\mu \left( x\right) =0
\label{front}
\end{equation}
and set 
\begin{equation}
h(w)=\frac{H\left( w\right) }{g(w)}  \label{h}
\end{equation}
where 
\begin{equation}
H\left( w\right) =\int_{0}^{1}y\left( x\right) f\left( wx\right) d\mu \left(
x\right) .  \label{hmais}
\end{equation}
The idea of the proof is to show that $h$ is constant and conclude from it
that $y$ must be null almost everywhere. The proof is not very long but, for
clarity purposes, we organize it in three straightforward steps.

\begin{description}
\item[Step 1] \emph{The function }$h$\emph{\ is entire and }$\varrho \left(
h\right) \leq 1.$
\end{description}

Because of its continuity, $f$ is bounded on every disk of the complex
plane. Therefore, the maximum of $f$ on a disk of radius $r$ exists and the
inequality 
\begin{equation}
M(r;H)\leq M(r;f)\left| \int_{0}^{1}y(x)d\mu \left( x\right) \right|
\label{inmax}
\end{equation}
holds. From this we infer that the integral defining $H$ converges uniformly
in compact sets. The condition (\ref{front}) forces every zero of $g$ to be
a zero of $H$ and the identity (\ref{h}) shows that $h$\ is an entire
function with less zeros than $H$; Since all functions are of order less
than one, then they can be written as canonical products. By (\ref{prop}),
the order of $h$ is less or equal to the order of $H$. On the other side,
the order of $H$ is less or equal to the order of $f$. This becomes clear
using (\ref{deford})\ and inequality (\ref{inmax}). It follows that $\rho
\left( h\right) \leq \rho (f)<1$.

\begin{description}
\item[Step 2] \emph{The function }$h$\emph{\ is constant.}
\end{description}

Condition (\ref{cond}) implies the existence of a constant $A>0$ such that $%
a_{n}\leq Ab_{n}$. Then $\left| x\right| \leq 1$ gives 
\begin{equation*}
f(itx)=\sum\limits_{n=0}^{\infty }a_{n}t^{2n}x^{2n}\leq
\sum\limits_{n=0}^{\infty }Ab_{n}t^{2n}x^{2n}\leq A\sum\limits_{n=0}^{\infty
}b_{n}t^{2n}=Ag(it).
\end{equation*}
Taking into account that $\mu $ is a positive-defined measure, this
inequality allows to estimate the integral in (\ref{h}) 
\begin{equation*}
\left| \int_{0}^{1}y(x)f(itx)d\mu (x)\right| \leq A\left| g(it)\right|
\int_{0}^{1}\left| y(x)\right| d\mu (x)
\end{equation*}
or equivalently 
\begin{equation*}
\left| h(it)\right| \leq A\int_{0}^{1}\left| y(x)\right| d\mu (x).
\end{equation*}
That is, $h$ is bounded on the imaginary axis. By Step1, $\varrho \left(
h\right) <1$. The Phragm\'{e}n-Lindel\"{o}f theorem with $\sigma =1$ shows
that $h$ is bounded in the complex plane. By Liouville theorem $h$ is a
constant.

\begin{description}
\item[Step 3] \emph{\ The function }$y\ $\emph{is null\ almost everywhere}$. 
$
\end{description}

Step 2 shows the existence of a constant $C$ such that $h\left( w\right) =C$
for every $w$ in the complex plane. Rewrite this as 
\begin{equation*}
\int_{0}^{1}y(x)f(wx)d\mu (x)-g\left( w\right) C=0
\end{equation*}
Use of the series expansion for $f(wx)$ and $g\left( w\right) $ gives 
\begin{equation*}
\sum\limits_{n=0}^{\infty }\left( -1\right) ^{n}\left[ a_{n}%
\int_{0}^{1}g(x)x^{2n}d\mu \left( x\right) -Cb_{n}\right] w^{2n}=0
\end{equation*}
by the identity theorem for analytical functions, 
\begin{equation}
\frac{a_{n}}{b_{n}}\int_{0}^{1}y(x)x^{2n}d\mu (x)=C.  \label{c1}
\end{equation}
On the other side, $x<1$ implies 
\begin{equation}
\left\vert \frac{a_{n}}{b_{n}}\int_{0}^{1}y(x)x^{2n}d\mu (x)\right\vert \leq 
\frac{a_{n}}{b_{n}}\left\vert \int_{0}^{1}y(x)d\mu (x)\right\vert .
\label{c2}
\end{equation}
Taking the limit when $n\rightarrow \infty $, (\ref{cond}) and (\ref{c1})
show that $C$ is null. As a result, for $n=1,2,...$, 
\begin{equation*}
\int_{0}^{1}y(x)x^{2n}d\mu (x)=0
\end{equation*}
Finally, the completeness of $x^{2n}$ in $L\left[ \mu ,\left( 0,1\right) %
\right] $ (by the M\"{u}ntz-Sz\'{a}sz theorem) shows that $y=0$ almost
everywhere. $\Box $

Before considering applications of theorem 1 it is convenient to recall that
if a function is given in its series form 
\begin{equation*}
f(z)=\sum\limits_{n=0}^{\infty }a_{n}z^{n\text{ }}
\end{equation*}%
then the order $\varrho \left( f\right) $ is given by 
\begin{equation}
\varrho \left( f\right) =\lim_{n\rightarrow \infty }\sup \frac{n\log n}{\log
\left( 1/\left\vert a_{n}\right\vert \right) }.  \label{ordlim}
\end{equation}

\subsection{$q$-special functions}

\subsubsection{Basic definitions}

Consider $0<q<1$. In what follows, the standard conventional notations from 
\cite{AAR} and \cite{GR}, will be used 
\begin{equation*}
(a;q)_{0}=1,\quad (a;q)_{n}=\prod_{k=1}^{n}(1-aq^{k-1}),
\end{equation*}%
\begin{equation*}
(a;q)_{\infty }=\lim\limits_{n\rightarrow \infty }(a;q)_{n},\quad
(a_{1},...,a_{m};q)_{n}=\prod_{l=1}^{m}(a_{l};q)_{n},\quad |q|<1,
\end{equation*}%
Jackson%
\'{}%
s $q$-integral in the interval $\left( 0,a\right) $ and in the interval $%
\left( 0,\infty \right) $ are defined, respectively, by 
\begin{eqnarray}
\int_{0}^{a}f\left( t\right) d_{q}t &=&\left( 1-q\right) a\sum_{n=0}^{\infty
}f\left( aq^{n}\right) q^{n}  \label{qinteg} \\
\int_{0}^{\infty }f\left( t\right) d_{q}t &=&\left( 1-q\right)
\sum_{n=-\infty }^{\infty }f\left( q^{n}\right) q^{n}\text{.}
\label{qinteginf}
\end{eqnarray}%
The $q$-difference operator $D_{q}$ is 
\begin{equation}
D_{q}f(x)=\frac{f(x)-f(qx)}{(1-q)x}.  \label{qdif}
\end{equation}%
These definitions appear in the formula of $q$- integration by parts 
\begin{equation}
\int_{0}^{1}G(qx)\left[ D_{q}f(x)\right] d_{q}x=f(1)G(1)-f(0)G(0)-%
\int_{0}^{1}f(x)D_{q}G(x)d_{q}x.  \label{qpartes}
\end{equation}%
We will denote by $L_{q}^{p}(X)$ the Banach space induced by the norm 
\begin{equation*}
\left\Vert f\right\Vert _{p}=\left[ \int_{X}\left\vert f\left( t\right)
\right\vert ^{p}d_{q}t\right] ^{\frac{1}{p}}.
\end{equation*}

There are three $q$-analogues of the Bessel function, all of them due to F.
H. Jackson and denoted by $J_{\nu }^{(1)}(x;q)$, $J_{\nu }^{(2)}(x;q)$ and $%
J_{\nu }^{(3)}(x;q)$. The third Jackson $q$-Bessel function has appeared
often in the literature under the heading \emph{The Hahn-Exton }$q$\emph{%
-Bessel function}. A well known formula usually credited to Hahn displays $%
J_{\nu }^{(2)}(x;q)$ as an analytical continuation of \ $J_{\nu }^{(1)}(x;q)$%
. Therefore, just the second and the third $q$-analogues are considered.
Their definition, in series form, is 
\begin{equation}
J_{\nu }^{\left( 2\right) }\left( x;q\right) =\frac{(q^{\nu +1};q)_{\infty }%
}{(q;q)_{\infty }}\sum_{n=0}^{\infty }\left( -1\right) ^{n}\frac{q^{n(\nu
+1)}}{(q^{\nu +1};q;q)_{n}}x^{2n+\nu }  \label{j2}
\end{equation}
\begin{equation}
J_{\nu }^{\left( 3\right) }\left( x;q\right) =\frac{(q^{\nu +1};q)_{\infty }%
}{(q;q)_{\infty }}\sum_{n=0}^{\infty }\left( -1\right) ^{n}\frac{q^{n(n+1)/2}%
}{(q^{\nu +1};q;q)_{n}}x^{2n+\nu }\text{.}  \label{j3}
\end{equation}
Very recently, Hayman \cite{H} obtained an asymptotic expansion for the
zeros of $J_{\nu }^{\left( 2\right) }$. For entire indices, the functions $%
J_{n}^{\left( 3\right) }(x;q)$ are generated by the relation, valid for $%
\left| xt\right| <1$, 
\begin{equation}
\frac{\left( qxt^{-1};q\right) _{\infty }}{\left( xt;q\right) _{\infty }}%
=\sum\limits_{n=-\infty }^{\infty }J_{n}^{\left( 3\right) }(x;q)t^{n}.
\label{gera}
\end{equation}
The Euler formula for the series form of an infinite product will be
critical on the remainder: 
\begin{equation}
\left( x;q\right) _{\infty }=\sum_{n=0}^{\infty }\left( -1\right) ^{n}\frac{%
q^{n\left( n-1\right) /2}}{(q;q)_{n}}x^{n}.  \label{euler}
\end{equation}

\subsubsection{Complete sets of $q$-special functions}

Theorem1 is very convenient to be applied to sets of $q$-special functions.
More often than not, these functions are of order zero, corresponding to the
situation where there is no restriction on the behavior of the zeros. The $q$%
-integral (\ref{qinteg}) is a Riemann-Stieltjes integral with respect to a
step function having infinitely many points of increase at the points $q^{k}$%
, with the jump at the point $q^{k}$ being $(1-q)q^{k}$.

Since (\ref{euler}) displays an easy relation between the zeros of a
function and its series form, we will use it first, for illustration
purposes, to construct complete (nonorthogonal) sets of infinite products.

\begin{example}
\textbf{\ }The sequence of infinite products $\left\{ (q^{-\frac{n}{2}%
+1}z^{2};q)_{\infty },n=0,1,...\right\} $\ forms a complete set in $%
L_{q}(0,1)$. For a proof of this take $f(z)=(qx^{2};q)_{\infty }$\ and $%
g(z)=(x^{2};q)_{\infty }$. Using Euler 
\'{}%
s Formula (\ref{euler}) one recognizes the setting of Theorem 1 with 
\begin{eqnarray*}
a_{n} &=&\frac{q^{\frac{n(n+1)}{2}}}{(q;q)_{n}} \\
b_{n} &=&\frac{q^{\frac{n(n-1)}{2}}}{(q;q)_{n}}.
\end{eqnarray*}%
Clearly, 
\begin{equation*}
\lim_{n\rightarrow \infty }\frac{a_{n}}{b_{n}}=0.
\end{equation*}%
Using (\ref{ordlim}) a short calculation recognizes $f$\ and $g$\ as
functions of order zero. By Theorem 1 it follows that $\{(q\lambda
_{n}z^{2};q)_{\infty }\}$\ is complete in $L_{q}(0,1)$, where $\lambda _{n}$%
\ is the $nth$\ zero of $(x^{2};q)_{\infty }$, that is, $\lambda _{n}=\pm
q^{-\frac{n}{2}}$, $n=0,1,..$.
\end{example}

Our next examples of complete sets, defined via the third Jackson $q$-Bessel
function, are of the form $J_{\nu }^{\left( 3\right) }\left( qx\lambda
_{n};q\right) $. It is well known that, denoting by $j_{n\nu }(q^{2})$ the $%
nth$ \ zero of $J_{\nu }^{\left( 3\right) }$, we have the orthogonality
relation 
\begin{equation}
\int_{0}^{1}xJ_{\nu }^{(3)}(qxj_{n\nu }(q^{2});q^{2})J_{\nu
}^{(3)}(qxj_{n\nu }(q^{2});q^{2})d_{q}x=0\text{,}  \label{ort}
\end{equation}
if $n\neq m$. It was proved in \cite{AB} that the system $\{x^{\frac{1}{2}%
}J_{\nu }^{(3)}(qxj_{n\nu };q^{2})\}$ forms an orthogonal basis of the space 
$L_{q}^{2}(0,1)$, by means of a $q$-version of a Dalzell criterion. In the
last section of \cite{AB} we proved the case (a) of the next theorem. The
case (b) will be used in section 4 of this paper.

\begin{example}
If $\nu >-1$, the sequence $\{J_{\nu }^{\left( 3\right) }\left( x\lambda
_{n};q^{2}\right) \}$\ is complete $L_{q}\left( 0,1\right) $\ if: (a) $%
\lambda _{n}=qj_{n,\alpha }^{\left( 3\right) }$, where $j_{n,\alpha
}^{\left( 3\right) }$\ is the $nth\ $zero of the function $J_{\alpha
}^{\left( 3\right) }\left( x;q^{2}\right) $; (b) $\alpha >-1$\ and $\lambda
_{n}=q^{-n}$, $n=0,1,...$Again we can build up the setting of theorem 1.
Indeed, it was proved in \cite{KS} that the roots of the third Jackson $q$%
-Bessel function are all real, simple and with countable cardinality. To
prove (a) consider $f$\ and $g$\ defined as 
\begin{eqnarray*}
\ f(x) &=&\frac{x^{-\nu }(q^{2};q^{2})_{\infty }}{(q^{2\nu
+2};q^{2})_{\infty }}J_{\nu }^{\left( 3\right) }\left( x;q^{2}\right) , \\
g(x) &=&\frac{x^{-\alpha }(q^{2};q^{2})_{\infty }}{(q^{2\alpha
+2};q^{2})_{\infty }}J_{\alpha }^{\left( 3\right) }(q^{-1}x;q^{2}).
\end{eqnarray*}%
Both $f$\ and $g$\ are functions of order $0$. Consequently, theorem 1 holds
with 
\begin{eqnarray*}
a_{n} &=&\frac{q^{n(n+1)}}{(q^{2u+2};q^{2};q^{2})_{n}}, \\
b_{n} &=&\frac{q^{n(n+1)-2n}}{(q^{2\alpha +2};q^{2};q^{2})_{n}}.
\end{eqnarray*}%
To prove the case (b) choose $f$\ as in (a) and $g(x)=(x^{2};q^{2})_{\infty
} $. Expand $g$\ by means of the series representation (\ref{euler}).\ The
result follows in a straightforward manner from theorem 1.
\end{example}

Now we will see the $J_{\nu }^{\left( 2\right) }\left( x;q\right) $ version
of the last example. The functions $J_{\nu }^{\left( 2\right) }\left(
x;q\right) $ are not orthogonal like the $J_{\nu }^{\left( 3\right) }\left(
x;q\right) $, but Rahman \cite{Rahman} was able to find a biorthogonality
relation, involving $J_{\nu }^{\left( 2\right) }\left( x;q\right) $ and $%
J_{\nu }^{\left( 1\right) }\left( x;q\right) $, that is reminiscent of (\ref%
{ort}). We will obtain the completeness property for the same range as in
Example 2. However, we will need a preliminary Lemma. The required Lemma is
the $q$-analogue of Theorem 5 in \cite{BP}.

\begin{lemma}
\ Let $\lambda _{n}$\ define a sequence of real numbers. For every $\nu >-1$%
, if the sequence $\{x^{-\nu -1}J_{\nu +1}^{(2)}(q\lambda _{n}x;q^{2})\}$\
is complete $L_{q}(0,1)$\ then $\{x^{-\nu }J_{\nu }^{(2)}(q\lambda
_{n}x;q^{2})\}$\ is also complete $L_{q}(0,1)$
\end{lemma}

\begin{proof}
Let $y(x)\in L_{q}(0,1)$ such that for every $n=1,2...$%
\begin{equation}
\int_{0}^{1}y(x)x^{-\nu }J_{\nu }^{(2)}(\lambda _{n}qx;q^{2})d_{q}x=0.
\label{c}
\end{equation}%
The $q$-difference operator (\ref{qdif}) acting on the power series (\ref{j2}%
) gives 
\begin{equation}
D_{q}\left[ x^{-\nu }J_{\nu }^{(2)}(\lambda _{n}x;q^{2})\right] =-\lambda
_{n}x^{-\nu }q^{\nu +1}J_{\nu +1}^{(2)}(\lambda _{n}xq;q^{2}).
\label{diffunc}
\end{equation}%
Now, use the $q$-integration by parts formula (\ref{qpartes})\ and (\ref%
{diffunc}) to obtain the identity 
\begin{eqnarray}
&&\int_{0}^{1}y(x)x^{-\nu }J_{\nu }^{(2)}(\lambda _{n}qx;q^{2})d_{q}x  \notag
\\
&=&q^{\nu +1}\lambda _{n}\int_{0}^{1}x^{-\nu -1}J_{\nu +1}^{(2)}(q\lambda
_{n}x;q^{2})\left[ x\int_{0}^{x}\left( q\lambda _{n}t\right) ^{\nu
}y(t)d_{q}t\right] d_{q}t\text{.}  \label{c4}
\end{eqnarray}%
By (\ref{c}), the expression (\ref{c4}) is zero for every $n=1,2,...$. Under
the hypothesis, $\{x^{-\nu -1}J_{\nu +1}^{(2)}(q\lambda _{n}x;q^{2})\}$ is
complete in $L_{q}(0,1)$. Clearly $x\int_{0}^{x}y(t)d_{q}t\in L_{q}(0,1)$
and thus, for $m=1,2,..,$ 
\begin{equation}
\int_{0}^{q^{m}}y(t)d_{q}t=0\text{.}  \label{lol}
\end{equation}%
This implies $y(q^{m})=0$ for every $m=1,2,..$.
\end{proof}

\begin{example}
\textbf{\ }\emph{I}f $\nu >-1$, the sequence $J_{\nu }^{\left( 2\right)
}\left( qx\lambda _{n};q^{2}\right) $\ is complete $L_{q}\left( 0,1\right) $%
\ if: $(a)\lambda _{n}=j_{n\alpha }^{\left( 2\right) },$ where $j_{n\alpha
}^{\left( 2\right) }\ $is the nth\ zero of the function $J_{\alpha }^{\left(
2\right) }\left( x;q^{2}\right) \ $and $\alpha >-1;$ (b) $\lambda
_{n}=q^{-n/2}$, $n=0,1,...$; First we remark that in \cite{I} the author
shows that the roots of the second Jackson $q$-Bessel function are all real
and simple and that there exists a countable infinite number of them. Then
use theorem 1 as in similar fashion as in the previous examples to establish
(a) when $\nu <\alpha +2$. A simple iteration of Lemma 1 yields the result
when $\alpha >-1$. On the other hand, (b) follows directly from theorem 1
choosing $g\left( z\right) =(z^{2};q^{4})_{\infty }.$
\end{example}

\section{Functions of order less than two}

Theorem 1 can be extended to the bigger class of entire functions of order
less than two. However, this requires a restriction on the behavior of the
zeros. With the same notational setting of the preceding theorem, the
following holds:

\begin{theorem}
If the order of $f$\ and $g$\ are less than two, then the sequence $%
\{f(\lambda _{n}x)\}$\ is complete $L\left[ \mu ,\left( 0,1\right) \right] $%
\ if, together with (\ref{cond}), the following condition holds\emph{\ } 
\begin{equation*}
\lambda _{n}\leq \zeta _{n}\text{.}
\end{equation*}
\end{theorem}

\begin{proof}
Consider $h$ defined as in (\ref{h}). The proof goes along the lines of the
proof of theorem 1. Only Step 2 requires a modification because now $\varrho
\left( h\right) <2$. The way to compensate this is to make the estimates
along smaller regions of the complex plane. Consider the angles defined by
the lines $\arg z=\pm \frac{\pi }{4}$ and $\arg z=\pm \frac{3\pi }{4}$.
These lines are the bounds of an angle of opening $\frac{\pi }{2}$. If $z$
belongs to one of the lines, then $z^{2}$ belongs to the imaginary axis. Say 
$z^{2}=it$, $t\in \mathbf{R}$. Now, by the Hadamard factorization theorem,
the infinite product expansion holds 
\begin{equation*}
\left\vert \frac{f(zx)}{g\left( z\right) }\right\vert =\prod_{n=1}^{\infty
}\left\vert \frac{\left( 1-\frac{itx^{2}}{\zeta _{n}^{2}}\right) }{\left( 1-%
\frac{it}{\lambda _{n}^{2}}\right) }\right\vert =\prod_{n=1}^{\infty }\left[ 
\frac{1+\frac{t^{2}x^{4}}{\zeta _{n}^{4}}}{1+\frac{t^{2}}{\lambda _{n}^{4}}}%
\right] ^{\frac{1}{2}}
\end{equation*}%
and the hypothesis $\lambda _{n}\leq \zeta _{n}$ together with $x\leq 1$
implies 
\begin{equation*}
\frac{1+\frac{t^{2}x^{4}}{\zeta _{n}^{4}}}{1+\frac{t^{2}}{\lambda _{n}^{4}}}%
\leq 1.
\end{equation*}%
Now, clearly 
\begin{equation*}
\left\vert f(zx)\right\vert \leq \left\vert g\left( z\right) \right\vert .
\end{equation*}%
From this we infer that the function $h$ is bounded on the sides of an angle
of opening $\frac{\pi }{2}$. Applying the Phragm\'{e}n-Lindel\"{o}f theorem
with $\sigma =2$ it follows that $h$ is bounded in the complex plane and, as
before, it is a constant.
\end{proof}

\subsection{Sets of Bessel functions}

Theorem 2 can be applied to the classical Bessel function. The Bessel
function of order $\nu >-1$ is defined by the power series 
\begin{equation*}
J_{\nu }(x)=\sum\limits_{n=0}^{\infty }\frac{(-1)^{n}}{n!\Gamma (\nu +n+1)}%
\left( \frac{x}{2}\right) ^{\nu +2n}\text{.}
\end{equation*}%
The function $(x/2)^{-\nu }J_{\nu }(x)$ is an entire function of order one
and it is well known that their zeros $\{j_{n\nu }\}$\ are all real and
simple.It is well known that the system $\{J_{\nu }(xj_{n\nu })\}$ is
orthogonal and complete in $L(0,1)$ and Boas and Pollard made in \cite{BP}
an extensive discussion of completeness properties of sets in the form $%
\{J_{\nu }(x\lambda _{n})\}$. We make yet another contribution to this topic
via the next example.

\begin{example}
Let $\alpha ,\nu >0$\ such that $\alpha $\ $<$\ $\nu $. The sequence $%
\left\{ J_{\nu }(xj_{n\alpha })\right\} $\ is then complete $L(0,1)$.
Consider $f(x)=(x/2)^{-\nu }J_{\nu }(x)$\ and $g(x)=\left( x/2\right)
^{-\alpha }J_{\alpha }(x)$. Both $f$\ and $g$\ are entire functions of the
form considered in Theorem 2, with 
\begin{eqnarray*}
a_{n} &=&\frac{1}{2^{2n}n!\Gamma (\nu +n+1)}, \\
b_{n} &=&\frac{1}{2^{2n}n!\Gamma (\alpha +n+1)}\text{.}
\end{eqnarray*}%
The identity $\Gamma (x+n+1)=\Gamma (x)(x)_{n+1}$\ implies 
\begin{equation*}
\frac{a_{n}}{b_{n}}=\frac{\Gamma (\alpha +n+1)}{\Gamma (\nu +n+1)}=\frac{%
\Gamma (\alpha )\left( \alpha \right) _{n+1}}{\Gamma (\nu )\left( \nu
\right) _{n+1}}\rightarrow 0\text{ }
\end{equation*}%
Furthermore, it is a well known fact from the theory of Bessel functions 
\cite[pag. 508]{Watson} that if $\alpha <\nu $\ then $j_{n\alpha }<j_{n\nu }$%
\ for all $n$.
\end{example}

\section{Uncertainty principles over $q$-linear grids}

As we have pointed out in the introduction, underlying the uncertainty
principle, there is the general idea that a function and its transform
cannot be both too small. A simple manifestation of this principle usually
occurs when a function $f$ and its transform $f\symbol{94}$ have both
bounded support (here we will consider the notion of support in an "almost
everywhere" sense: a function is said to be supported on a set $A$ if it
vanishes almost everywhere outside $A$). If the measure is continuous, then
the transform is analytic and vanishes in a set with an accumulation point.
Therefore it must vanish identically. If we have an inversion formula then
also $f$ vanishes identically. This is the case of the Fourier and the
Hankel transform.

When dealing with discrete versions of uncertainty principles one often
fiunds changes that go beyond formal considerations. For instance, since we
are dealing with almost everywhere supports, the discrete analogue of
"vanishing outside an interval" is "vanishing at the points that support a
discrete measure outside an interval". The analytic function argument used
above will then fail, given the measure has no accumulation point outside
the interval. It may simply happen that the transform is not vanishing in a
sufficiently coarse set (in particular, no accumulation point) to make the
function vanish.

A particulary significative example occurs when considering a measure
supported on the integer powers of a real number $q\in \lbrack 0,1]$, like
Jackson 
\'{}%
s $q$-integral between $0$ and $\infty $. The support of the measure is $%
\left\{ q^{n}\right\} _{n=-\infty }^{\infty }$ which has zero as the only
accumulation point. Split this support in two grids: The one with positive
powers of $q$ accumulates at zero. The other consists in negative powers of $%
q$, the gap between the points increasing at a geometrical rate. Given the
sparsity of the grid $\{q^{n}\}_{n=-\infty }^{0}$, we might be skeptical
about the fact that simultaneous vanishing of $f$ and its $q$-discrete
transform in such a set is enough to force vanishing at the remaining
support points of the measure. However, in the case of a certain $q$%
-discrete transform whose kernel is the third Jackson $q$-Bessel function,
we will see in theorem 3 that this is indeed the case.

\subsection{The $q$-Hankel transform}

Follow Koornwinder and Swarttouw \cite{KS}, and define a $q$-Hankel
transform setting 
\begin{equation}
(H_{q}^{\nu }f)\left( x\right) =\int_{0}^{\infty }\left( xt\right) ^{\frac{1%
}{2}}J_{\nu }^{(3)}(xt;q^{2})f\left( t\right) d_{q}t\text{.}  \label{qhank}
\end{equation}%
It was shown in \cite{KS} that the $q$-Hankel transform satisfies the
inversion formula\emph{\ } 
\begin{equation}
f\left( t\right) =\int_{0}^{\infty }\left( xt\right) ^{\frac{1}{2}%
}(H_{q}^{\nu }f)\left( x\right) J_{\nu }^{(3)}(xt;q^{2})d_{q}x=(H_{q}^{\nu
}(H_{q}^{\nu }f))\left( t\right)  \label{hankelinv}
\end{equation}%
where $t$\ takes the values\emph{\ }$q^{k},k=0,\pm 1,\pm 2,...$\emph{. }%
Since the transform $H_{q}^{\nu }$ is self-inverse, it provides a Hilbert
space isometry between $L_{q}^{2}(0,1)$ and the space\ 
\begin{equation}
PW_{q}^{\nu }=\left\{ f\in L_{q}^{2}(\mathbf{R}^{+}):f\left( x\right)
=\int_{0}^{1}\left( tx\right) ^{\frac{1}{2}}J_{\nu }^{(3)}(xt;q^{2})u\left(
t\right) d_{q}t,u\in L_{q}^{2}\left( 0,1\right) \right\} \text{.}
\label{qespaco}
\end{equation}%
This space was defined in \cite{A}\ as the $q$-Bessel version of the
Paley-Wiener space of bandlimited functions and it was recognized as being a
reproducing kernel Hilbert space, with an associated $q$-sampling theorem.

\subsection{A vanishing theorem for the $q$-Hankel transform}

The vanishing theorem for the $q$-Hankel transform is now a simple
consequence of the completeness result on sets of third $q$-Bessel functions.

\begin{theorem}
Let $f\in L_{q}(R^{+})$\ such that both $f$\ and its $q$-Hankel transform
vanish at the points $q^{-n},n=0,1,...$. Then 
\begin{equation*}
f(q^{k})=0,k=0,\pm 1,\pm 2,...
\end{equation*}%
that is, $f$\ vanishes in the equivalent classes of $L_{q}(R^{+})$.
\end{theorem}

\begin{proof}
\textbf{\ }Let $f\in L_{q}(\mathbf{R}^{+})$. If $f(q^{-n})=0,n=0,1,...$,
then the $q$-Hankel transform of $f$\ is 
\begin{equation}
H_{q}^{\nu }f(\omega )=\int_{0}^{1}\left( \omega t\right) ^{\frac{1}{2}%
}J_{\nu }^{\left( 3\right) }(\omega t;q^{2})f\left( t\right) d_{q}t\text{.}
\label{h01}
\end{equation}%
Since our second assumption says that $(H_{q}^{\nu }f)(q^{-n})=0,n=0,1,...$,
if we set $\ \omega =q^{-n}$ in (\ref{h01}), the result is 
\begin{equation}
\int_{0}^{1}(q^{-n}t)^{\frac{1}{2}}J_{\nu }^{\left( 3\right)
}(q^{-n}t;q^{2})f(t)d_{q}t=0,n=0,1,...\text{.}  \label{h02}
\end{equation}%
Now we have from Example 3 (b) that, if $\nu >-1$, then the sequence $%
\{J_{\nu }^{\left( 3\right) }\left( q^{-n}t;q^{2}\right) \}$ is complete in $%
L_{q}^{1}\left( 0,1\right) $. Therefore, the conditions (\ref{h02}) imply $%
f\equiv 0$ in $L_{q}^{1}\left( 0,1\right) $, that is, $f(q^{n})=0,n=0,1,...$%
. Since, by hypothesis, $f(q^{-n})=0,n=0,1,...$, the result follows.
\end{proof}

The vanishing theorem has a prompt consequence when seen in terms of $%
PW_{q}^{\nu }$.

\begin{corollary}
\ $\Gamma =\{q^{-n},n\in N\}$\ is a set of uniqueness for the space $%
PW_{q}^{\nu }$.
\end{corollary}

\begin{proof}
Take $f\in PW_{q}^{\nu }$ such that $f(q^{-n})=0,n=1,2,...$. If $f$ is of
the form required in (\ref{qespaco}) then $f=H_{q}^{\nu }u^{\ast }$ where $%
u^{\ast }\in L_{q}^{2}\left( \mathbf{R}^{+}\right) $ is obtained from $u\in
L_{q}^{2}\left( 0,1\right) $ by prescribing $u(q^{-n})=0,n\in \mathbf{N}$.
By the inversion formula (\ref{hankelinv}), $u^{\ast }=H_{q}^{\nu }f$. We
conclude that $H_{q}^{\nu }f(q^{-n})=0,n=0,1,...$. By Theorem 3, $f\equiv 0.$
\end{proof}

\begin{remark}
Observe that, if $\left( H_{q}^{\nu }f\right) \left( q^{-n}\right) =0,n\in N$%
, taking into account definitions (\ref{qinteg})\ and (\ref{qinteginf}) then 
$f=\left( H_{q}^{\nu }\left( H_{q}^{\nu }f\right) \right) $\ is of the form
required in (\ref{qespaco}). The argument in proof of Corollary 1 shows the
following characterization of $PW_{q}^{\nu }$: 
\begin{equation*}
PW_{q}^{\nu }=\left\{ f\in L_{q}^{2}(\mathbf{R}^{+}):(H_{q}^{\nu
}f)(q^{-n})=0,n=1,2,...\right\} \text{.}
\end{equation*}%
The property $\left( H_{q}^{\nu }f\right) \left( q^{-n}\right) =0,n=0,1,...$%
\ can thus be seen as a sort of \textquotedblright $q$%
-Hankel-bandlimitedness\textquotedblright . It was shown in \cite{A} that
there are many features in this space analogous to the classical Paley
Wiener space, including a sampling theorem and a reproducing kernel.\ 
\end{remark}

\begin{remark}
If $\nu >0$, $y>-\frac{1}{2}$\ and $x\in R$, the following $q$-analogue of
the Sonine integral was proved in \cite{A}: 
\begin{equation}
\frac{(q;q)_{\infty }}{(q^{\nu };q)_{\infty }}x^{-\nu }J_{y+\nu
}(x;q)=\int_{0}^{1}t^{\frac{y}{2}}\frac{(tq;q)_{\infty }}{(tq^{\nu
};q)_{\infty }}J_{y}(xt^{\frac{1}{2}};q)d_{q}t\text{.}  \label{qsonine}
\end{equation}%
using this formula we have seen that, if $\alpha >\nu >-\frac{1}{2}$, the
function $f(x)=x^{\nu -\alpha +\frac{1}{2}}J_{\alpha }(x;q^{2})$\ belongs to
the space $PW_{q}^{\nu }$\ and its image via the $q$-Hankel transform, in
the space $L_{q}^{2}(0,1)$, is the function 
\begin{equation*}
u(t)=(1+q)t^{\nu +\frac{1}{2}}\frac{(q^{2\alpha -2\nu };q^{2})_{\infty
}(t^{2}q^{2};q^{2})_{\infty }}{(q^{2};q^{2})_{\infty }(t^{2}q^{2\alpha -2\nu
};q^{2})_{\infty }}.
\end{equation*}%
Observe that the condition $\left( H_{q}^{\nu }f\right) \left( q^{-n}\right)
=u\left( q^{-n}\right) =0,n=1,2,...$\ is clearly satisfied. The function $f$%
\ is thus an example of a $q$-bandlimited function.
\end{remark}

\begin{remark}
A signal theoretical interpretation follows from Corollary 1: If we identify
a function $f$\ with its representant in the equivalent classes of $%
L^{1}(R^{+})$, we can think about $f$\ as a discrete signal with points 
\begin{equation*}
\{...,f(q^{n}),...,f(q),f(1),f(q^{-1}),...,f(q^{-n}),...\}.
\end{equation*}%
If a signal $f\in PW_{q}^{\nu }$\ is transmitted along a channel, and any
set of points contained in $\{...,f(q^{n}),...,f(q),f(1)\}$\ is
\textquotedblright lost\textquotedblright , then the received signal $g$\
still contains the whole information about $f$. To prove this, observe that $%
f-g\in PW_{q}^{\nu }$\ vanishes on a set containing $\Gamma =\{q^{-n},n\in
N\}$. By Corollary 1, $f=g$. Similar ideas were explored in \cite{DS}, in
the context of bandlimited signals with missing segments on time domain and
in the recovery of sparse discrete finite signals with missing samples.
\end{remark}

\subsection{An uncertainty principle with $\protect\epsilon $ \ 
\'{}%
s}

To complete our discussion on uncertainty principles for the $q$-Hankel
transform, we will now lose contact with the completeness concept that has
been our unifying theme until so far, and borrow ideas from modern signal
analysis.

The notion of $\epsilon $-concentration is required in order to obtain
information of a more quantitative character. A function $f\in L^{2}(X,\mu )$
such that $\left\| f\right\| _{L^{2}(X,\mu )}=1$ is said to be $\epsilon
_{T} $-concentrated in a set $T$ if 
\begin{equation}
\left\| f-\mathbf{1}_{T}\right\| _{L^{2}(X,\mu )}\leq \epsilon _{T}\text{.}
\label{epsilon}
\end{equation}
In \cite{DS}, Donoho and Stark proved that if a function $f$ of unit $L^{2}(%
\mathbf{R})$ norm is $\epsilon _{T}$-concentrated in a measurable set $T$
and its Fourier transform is $\epsilon _{\Omega }$-concentrated in a
measurable set $\Omega $, then $\left| T\right| \left| \Omega \right| \geq
(1-\epsilon _{T}-\epsilon _{\Omega })^{2}$, where $\left| {}\right| $
denotes Lebesgue measure. The uncertainty principle of Donoho and Stark was
extended by de Jeu to general bounded integral operators satisfying a
Plancherel theorem \textbf{\cite{Jeu}}. De Jeu 
\'{}%
s result is of a very general scope, and it will be stated here in the
degree of generality suitable to our needs.

\textbf{Theorem A \cite{Jeu}} \emph{Consider an integral transform defined,
for every }$f\in L^{2}(X,\mu )$ \emph{by }$(Kf)(x)=\int_{X}K(x,t)f(t)d\mu
(x) $\emph{,} \emph{mapping }$L^{2}(X,\mu )$ \emph{in itself, and such that
there is a Plancherel theorem for all its range. If }$f$ \emph{is of unit
norm and }$\epsilon _{T}$\emph{-concentrated in} $T$\emph{\ and }$Kf$ \emph{%
is} $\epsilon _{\Omega }$\emph{-concentrated in} $\Omega $, \emph{then the
following inequality holds:} 
\begin{equation}
\left\Vert \mathbf{1}_{T\times \Omega }K(x,t)\right\Vert _{L^{2}(X,\mu
)}\geq 1-\epsilon _{T}-\epsilon _{\Omega }\text{.}  \label{in0}
\end{equation}

It is possible to use Theorem A\textbf{\ }to extract more valuable
information about the size of the $\epsilon $-concentration sets in the case
of the $q$-Hankel transform and obtain an uncertainty principle of Donoho
and Stark style for $\epsilon $-concentration in sets of the form $%
T=\{q^{n+n_{T}}\}_{n=0}^{\infty }$.

Observe that the $q$-integral over the set $\left( 0,q^{n_{T}}\right) $, $%
n_{T}\in \mathbf{Z}$, is 
\begin{equation}
\int_{0}^{q^{n_{T}}}f\left( t\right) d_{q}t=\left( 1-q\right)
\sum_{n=0}^{\infty }f\left( q^{n+n_{T}}\right) q^{n+n_{T}}\text{.}
\label{qint}
\end{equation}
and $\epsilon _{T}$-concentration in a set $T=\{q^{n+n_{T}}\}_{n=0}^{\infty
} $\ in the $L_{q}^{2}(\mathbf{R}^{+})$ norm becomes, attending to (\ref%
{qinteginf}) and (\ref{epsilon}), 
\begin{equation*}
\left( 1-q\right) ^{\frac{1}{2}}\sum_{n=-\infty }^{n_{A}+1}\left\vert
f\left( q^{n}\right) \right\vert ^{2}q^{2n}\leq \epsilon _{T}^{2}\text{.}
\end{equation*}

The uncertainty principle in this context reads as follows.

\begin{theorem}
\emph{If }$f\in L_{q}^{2}(R^{+})$\emph{\ of unit norm is }$\epsilon _{T}$%
\emph{-concentrated in }$\{q^{n+n_{T}}\}_{n=0}^{\infty }$\emph{\ and }$%
H_{q}^{\nu }f$\emph{\ is }$\epsilon _{\Omega }$\emph{-concentrated in }$%
\{q^{n+n_{\Omega }}\}_{n=0}^{\infty }$\emph{, then } 
\begin{equation*}
n_{T}+n_{\Omega }\geq 2\log _{q}[(q^{2};q^{2})_{\infty }^{2}(1-\epsilon
_{T}-\epsilon _{\Omega })]\text{.}
\end{equation*}
\end{theorem}

\begin{proof}
\textbf{\ }The proof uses the following estimate of the third Jackson $q$%
-Bessel function obtained in \cite{FHB}. For every $x=q^{k}$, $k=0,\pm
1,2,...$ the inequality holds 
\begin{equation}
\left\vert J_{\nu }^{(3)}(x;q)\right\vert \leq \frac{x^{\nu }}{%
(q;q^{2})_{\infty }^{2}}\text{.}  \label{in2}
\end{equation}%
Now observe that if $q^{n_{T}+n_{\Omega }}\geq 1$ then the proposition is
trivial, since $(q^{2};q^{2})_{\infty }^{2}<1$. Thus we can assume without
loss of generalization that $q^{n_{T}+n_{\Omega }}<1$. In this case we have $%
xt=q^{k}$ for some entire $k$. Then, use of (\ref{in2})\ together with the
definition of the $q$-integral yields, after applying Theorem A to the $q$%
-Hankel transform gives 
\begin{eqnarray*}
1-\epsilon _{T}-\epsilon _{\Omega } &\leq &\left\Vert \mathbf{1}_{\left[
0,q^{n_{T}}\right] \times \left[ 0,q^{n_{\Omega }}\right] }(x,t)\left(
xt\right) ^{\frac{1}{2}}J_{\nu }^{(3)}\left( xt;q^{2}\right) \right\Vert
_{L_{q}^{2}(\mathbf{R}^{+})\times L_{q}^{2}(\mathbf{R}^{+})} \\
&=&\int_{0}^{q^{n_{\Omega }}}\left[ \int_{0}^{q^{n_{T}}}\left[ \left(
tx\right) ^{\frac{1}{2}}J_{\nu }^{(3)}\left( xt;q^{2}\right) \right]
^{2}d_{q}t\right] d_{q}x \\
&\leq &\int_{0}^{q^{n_{\Omega }}}\int_{0}^{q^{n_{T}}}\left[ \frac{1}{%
(q;q^{2})_{\infty }^{2}}\right] ^{2}d_{q}td_{q}x \\
&=&\frac{q^{n_{T}+n_{\Omega }}}{(q;q^{2})_{\infty }^{2}}\text{.}
\end{eqnarray*}
\end{proof}

\textbf{Acknowledgement. }\emph{This work was partially supported by CMUC
and FCT post-doctoral grant SFRH/BPD/26078/2005. Some of the ideas in the
first part were included in the work I did under supervision of Joaquin
Bustoz at Arizona State University. The last part of it was done at NuHAG,
University of Vienna, hosted by Hans Feichtinger. }

\end{document}